\ifx\documentclass\undefined
\documentstyle[12pt]{article}
\else
\documentclass[12pt]{article}
\usepackage{latexsym}
\usepackage{amsfonts}
\usepackage{amsmath}
\usepackage{amssymb}
\fi

\author{ K\'aroly Bezdek 
\thanks{Partially supported by a Natural Sciences and 
Engineering Research Council of Canada Discovery Grant and by the Hung. Acad. Sci. Found. (OTKA), grant no. T72537. (This survey is partially based on the author's talk delivered at the meeting "Intuitive Geometry, in Memoriam L\'aszl\'o Fejes T\'oth", June 30 - July 4, 2008, Budapest, Hungary.)}}

\font\tenBbb=msbm10 at 12pt         \font\sevenBbb=msbm9    \font\fiveBbb=msbm7
\newfam\Bbbfam
\textfont\Bbbfam=\tenBbb \scriptfont\Bbbfam=\sevenBbb
\scriptscriptfont\Bbbfam=\fiveBbb

\def\E{{\mathbb E}}
\def\H{{\mathbb H}}
\def\S{{\mathbb S}}

\def\kkk{\null\hfill $\Box$\smallskip}
\def\r{ \right}

\def\la{\left\langle}
\def\ra{\r\rangle}

\def\KK{{\bf K}}
\def\LL{{\bf L}}

\def\PP{{\bf P}}
\def\BB{{\bf B}}
\def\CC{{\bf C}}

\def\QQ{{\bf Q}}
\def\TT{{\bf T}}

\newcommand{\proof}{{\noindent\bf Proof:{\ \ }}}

\newtheorem{theorem}{Theorem}[section]
\newtheorem{lemma}[theorem]{Lemma}

\newtheorem{remark}[theorem]{Remark}

\newtheorem{cor}[theorem]{Corollary}
\newtheorem{con}[theorem]{Conjecture}
\newtheorem{prob}[theorem]{Problem}

\newtheorem{definition}[theorem]{Definition}

\newcommand{\vol}{\mathop{\rm vol\,}}

\newcommand{\crv}{{\rm crv}}

\title{Tarski's plank problem revisited
}

\begin{document}

\maketitle

\date

\begin{abstract}
In the 1930's, Tarski introduced his plank problem at a time when the field Discrete Geometry was about to born. It is quite remarkable that Tarski's question and its variants continue to generate interest in the geometric and analytic aspects of coverings by planks in the present time as well. The paper is a survey type with a list of open research problems.
\end{abstract}

\medskip

\section{Introduction}

Tarski's plank problem has generated a great interest in understanding the geometry of coverings by planks. There have been a good number of results published in connection with the plank problem of Tarski that are surveyed in this paper. The paper is divided into six sections entitled Plank theorems; Covering convex bodies by cylinders; Covering lattice points by hyperplanes; Partial coverings by planks; Linking the Kakeya-P\'al and the Blaschke-Lebesgue problems to the partial covering problem and Strengthening the plank theorems of Ball and Bang. Also, there are a number of open problems mentioned in order to motivate further research.  

\medskip

\section{Plank theorems}

A {\it convex body} of the Euclidean space $\E^d$ is a compact convex set with non-empty interior. Let $\CC\subset\E^d$ be a convex body, and let $H\subset\E^d$ be a hyperplane. Then the distance ${\rm w}(\CC , H)$ between the two supporting hyperplanes of $\CC$ parallel to $H$ is called the {\it width of $\CC$ parallel to $H$}. Moreover, the smallest width of $\CC$ is called the {\it minimal width} of $\CC$ and is denoted by ${\rm w}(\CC)$.  

Recall that in the 1930's, Tarski posed what came to be known as the plank problem. A {\it plank} $\PP$ in $\E^d$ is the (closed) set of points between two distinct parallel hyperplanes. The {\it width} ${\rm w}(\PP)$ of $\PP$ is simply the distance between the two boundary hyperplanes of $\PP$. Tarski conjectured that if a convex body of minimal width $w$ is covered by a collection of planks in $\E^d$, then the sum of the widths of these planks is at least $w$. This conjecture was proved by Bang in his memorable paper \cite{Ba51}. (In fact, the proof presented in that paper is a simplification and generalization of the proof published by Bang  somewhat earlier in \cite{Ba50}.) Thus, the following statement we call the plank theorem of Bang.

\begin{theorem}
If the convex body $\CC$ is covered by the planks $\PP_1, \PP_2, \dots , \PP_n$ in $\E^d$ (i.e. $\CC\subset \PP_1\cup \PP_2\cup \dots \cup \PP_n\subset\E^d$), then
$$\sum_{i=1}^n {\rm w}(\PP_i)\ge {\rm w}(\CC).$$
\end{theorem}

In \cite{Ba51}, Bang has raised the following stronger version of Tarski's plank problem. We phrase it via the following definition.
Let $\CC$ be a convex body and let $\PP$ be a plank with boundary hyperplanes parallel to the hyperplane $H$ in $\E^d$. We define the {\it $\CC$-width} of the plank $\PP$ as $\frac{{\rm w}(\PP) }{{\rm w}(\CC , H) }$ and label it by ${\rm w}_{\CC}(\PP)$. 

\begin{con}\label{Bang-conjecture}
If the convex body $\CC$ is covered by the planks $\PP_1, \PP_2, \dots , \PP_n$ in $\E^d$, then 
$$\sum_{i=1}^n {\rm w}_{\CC}(\PP_i)\ge 1.$$
\end{con}

It was Alexander \cite{Al68}, who noticed, that Conjecture~\ref{Bang-conjecture} is equivalent to the following
generalization of a problem of Davenport.

\begin{con}\label{Alexander-Davenport}
If a convex body $\CC$ in $\E^d$ is sliced by $n-1$ hyperplane cuts (pieces are held together after each cut), then there exists a piece that covers a translate of $\frac{1}{n}\CC$.
\end{con}

We note that that in \cite{BeBe96} one can find a proof of  Conjecture~\ref{Alexander-Davenport} for successive hyperplane cuts (i.e. for hyperplane cuts when one piece is divided by each cut). Also, the same paper (\cite{BeBe96}) introduced two additional equivalent versions for Conjecture~\ref{Bang-conjecture}. As they seem to be of independent interest we recall them following the terminology used in \cite{BeBe96}.

Let $\CC$ and $\KK$ be convex bodies in $\E^d$ and let $H$ be a hyperplane of $\E^d$. The {\it $\CC$-width of $\KK$ parallel to $H$} is denoted by $ {\rm w}_{\CC}(\KK , H)$ and is defined as $\frac{{\rm w}(\KK , H)}{{\rm w}(\CC , H) }$. The {\it minimal $\CC$-width of $\KK$} is denoted by $ {\rm w}_{\CC}(\KK)$ and is defined as the minimum of $ {\rm w}_{\CC}(\KK , H)$, where the minimum is taken over all possible hyperplanes $H$ of $\E^d$. Recall that the inradius of $\KK$ is the radius of the largest ball contained in $\KK$. It is quite natural then to introduce the {\it $\CC$-inradius of $\KK$} as the factor of the largest (positively) homothetic copy of $\CC$, a translate of which is contained in $\KK$. We need to do one more step to introduce the so-called successive $\CC$-inradii of $\KK$ as follows. Let $r$ be the $\CC$-inradius of $\KK$. For any $0<\rho\le r$ let the {\it $\rho\CC$-rounded body of $\KK$} be denoted by ${\KK}^{\rho\CC}$ and be defined as the union of all translates of $\rho\CC$ that are covered by $\KK$. Now, take a fixed integer $n\ge 1$. On the one hand, if $\rho>0$ is sufficiently small, then $ {\rm w}_{\CC}({\KK}^{\rho\CC})>n\rho$. On the other hand, $ {\rm w}_{\CC}({\KK}^{r\CC})=r\le nr$. As $ {\rm w}_{\CC}({\KK}^{\rho\CC})$ is a decreasing continuous function of $\rho>0$ and $n\rho$ is a strictly increasing continuous function of $\rho$ there exists a uniquely determined $\rho>0$ such that 
$$ {\rm w}_{\CC}({\KK}^{\rho\CC})=n\rho.$$

This uniquely determined $\rho$ is called the {\it $n$-th successive $\CC$-inradius of $\KK$} and is denoted by $r_{\CC}(\KK , n)$. Notice that $r_{\CC}(\KK , 1)=r$. Now, the two equivalent versions of Conjecture~\ref{Bang-conjecture} and Conjecture~\ref{Alexander-Davenport} introduced in \cite{BeBe96} can be phrased as follows.

\begin{con}\label{Bezdek-Bezdek-1}
If a convex body $\KK$ in $\E^d$ is covered by the planks $\PP_1,$ $\PP_2, \dots , \PP_n$, then $\sum_{i=1}^n {\rm w}_{\CC}(\PP_i)\ge
{\rm w}_{\CC}(\KK)$ for any convex body $\CC$ in $\E^d$.
\end{con}

\begin{con}\label{Bezdek-Bezdek-2}
Let $\KK$ and $\CC$ be convex bodies in $\E^d$. If $\KK$ is sliced by $n-1$ hyperplanes (pieces are held together after each cut), then the minimum of the greatest $\CC$-inradius of the pieces is equal to the $n$-th successive $\CC$-inradius of $\KK$, i.e. it is $r_{\CC}(\KK , n)$.
\end{con}

The main result of \cite{BeBe96} is the next theorem that (under the condition that $\CC$ is a ball) answers a question raised by Conway (\cite{BeBe95}) as well as proves Conjecture~\ref{Bezdek-Bezdek-2} for successive hyperplane cuts.

\begin{theorem}
Let $\KK$ and $\CC$ be convex bodies in $\E^d, d\ge 2$. If $\KK$ is sliced into $n$ pieces by $n-1$ successive hyperplane cuts (i.e. when one piece is divided by each cut), then the minimum of the greatest $\CC$-inradius of the pieces is the $n$-th successive $\CC$-inradius of $\KK$, i.e. $r_{\CC}(\KK , n)$. An optimal partition is achieved by $n-1$ parallel hyperplane cuts equally spaced along the $\CC$-width of the $r_{\CC}(\KK , n)\CC$-rounded body of $\KK$.
\end{theorem}

The special case of Conjectures~\ref{Bang-conjecture}, ~\ref{Alexander-Davenport}, ~\ref{Bezdek-Bezdek-1}, ~\ref{Bezdek-Bezdek-2}, when the convex body to be covered or to be sliced is centrally symmetric, has been proved by Ball in his celebrated paper \cite{B91}. Thus, the following is the plank theorem of Ball.

\begin{theorem}
If a centrally symmetric convex body $\CC$ in $\E^d$ is sliced by $n-1$ hyperplane cuts (pieces are held together after each cut), then there exists a piece that covers a translate of $\frac{1}{n}\CC$.
\end{theorem}

We close this section by mentioning that Ball \cite{B01} has dramatically improved his plank theorem for complex Hilbert spaces. 

\medskip

\section{Covering convex bodies by cylinders}

\medskip

The main problem of this section is due Bang(1951). In his paper \cite{Ba51}, Bang raises the following challanging question.

\begin{prob}\label{BBang}
Prove or disprove that the sum of the base areas of finitely many cylinders covering a $3$-dimensional convex body is at least half of the minimum area 2-dimensional projection of the body.
\end{prob}

If true, then Bang's estimate is sharp due to a covering of a regular tetrahedron by two cylinders described in \cite{Ba51}. A very recent paper of Litvak and the author (\cite{BeLi09})
investigates Problem~\ref{BBang} as well as its higher dimensional analogue. Their main result can be summarized as follows. 

Given $0< k<d$ define a $k$-codimensional cylinder $\CC$ in $\E^d$ as a set which can be presented in the 
form $\CC = H + B$, where $H$ is a $k$-dimensional linear subspace of $\E^d$ and  $B$ is a measurable set in the orthogonal complement $H^{\perp}$ of $H$.  For a given convex body $\KK$ 
and a $k$-codimensional cylinder $\CC= H + B$ we define the cross-sectional volume $\crv _{\KK} (\CC)$ of $\CC$ with respect to $\KK$ as follows 
$$
 \crv _{\KK} (\CC) := \frac{\vol _{d-k} (\CC\cap H^{\perp})}{\vol _{d-k} (P_{H^{\perp}} 
 \KK)} =\frac{\vol _{d-k} (P_{H^{\perp}} \CC)}{\vol _{d-k} (P_{H^{\perp}}
 \KK)} = \frac{\vol _{d-k} (B)}{\vol _{d-k} (P_{H^{\perp}} \KK)} ,
$$
where $P_{H^{\perp}}: \E^d \to H^{\perp}$ denotes the orthogonal projection of $\E^d$ onto $H^{\perp}$. Notice that for every invertible affine map $T: \E^d \to \E^d$ one has $\crv _{\KK} (\CC) = \crv _{T \KK} (T \CC)$. The following theorem is proved 
in \cite{BeLi09}.

\begin{theorem}\label{Bezdek-Litvak}
Let $\KK$ be a convex body in $\E^d$. Let $\CC_1, \dots, \CC_N$ 
be $k$-co\-di\-men\-si\-o\-nal cylinders in $\E^d, 0< k<d$ such that 
$\KK\subset \bigcup _{i=1}^N \CC_i .$ Then 
$$
    \sum _{i=1}^N  \mbox{\rm crv} _{\KK} (\CC_i) \geq  \frac{1}{{d \choose k}}.
$$ 
Moreover, if $\KK$ is an ellipsoid and $\CC_1, \dots, \CC_N$ 
are $1$-co\-di\-men\-si\-o\-nal cylinders in $\E^d$ such that $\KK\subset \bigcup _{i=1}^N \CC_i$, then 
$$
    \sum _{i=1}^N  \mbox{\rm crv} _{\KK} (\CC_i) \geq 1.  
$$
\end{theorem}

The case $k=d-1$ of Theorem~\ref{Bezdek-Litvak} corresponds to Bang's affine invariant version of Tarski's plank problem (\cite{Ba51}). Indeed, in this case we have the sum of relative widths of the body. Recall that Ball (\cite{B91}) proved that such sum should exceed $1$ in the case of centrally symmetric body $\KK$, while the general case is still open. Theorem~\ref{Bezdek-Litvak} for $k=d-1$ implies the lower bound $1/d$. 

As an immediate corollary of Theorem~\ref{Bezdek-Litvak} we get the following estimate for Problem~\ref{BBang}.

\begin{cor}
The sum of the base areas of finitely many ($1$-codimensional) cylinders covering a $3$-dimensional convex body is always at least one third of the minimum area $2$-dimensional projection of the body.
\end{cor}

Also, note that the inequality of Theorem~\ref{Bezdek-Litvak} on covering ellipsoids by $1$-codimensional cylinders is best possible. By looking at this result from the point of view of $k$-codimensional cylinders we are led to ask the following still open question.

\begin{prob}
If $\KK$ is an ellipsoid and $\CC_1, \dots, \CC_N$ 
are $k$-co\-di\-men\-si\-o\-nal cylinders in $\E^d, 1<k<d-1$ such that $\KK\subset \bigcup _{i=1}^N \CC_i$, then prove or disprove that 
$$
    \sum _{i=1}^N  \mbox{\rm crv} _{\KK} (\CC_i) \geq 1.  
$$
\end{prob}

\section{Covering lattice points by hyperplanes}

\medskip

In their paper \cite{BeHa91}, Hausel and the author have established the following discrete version of Tarski's plank problem. 

Recall that the lattice width of a convex body $\KK$ in $\E^d$ is defined as
$$
  w(\KK, \mathbb Z ^d ) = \min\left\{ \max_{x\in \KK} \la x, y \ra - 
  \min_{x\in \KK} \la x, y \ra  \ \mid \ y\in \mathbb Z ^d, 
  \ y\ne 0 \r\},
$$
where $\mathbb Z ^d$ denotes the integer latice of $\E^d$. It is well-known that if $y\in \mathbb Z ^d, 
  \ y\ne 0$ is chosen such that $\lambda y\notin \mathbb Z ^d$ for any $0<\lambda <1$, then 
$$ \max_{x\in \KK} \la x, y \ra - \min_{x\in \KK} \la x, y \ra$$ 
is equal to the Euclidean width of $\KK$ in the direction $y$ divided by the Euclidean distance between two consecutive lattice hyperplanes of $\mathbb Z ^d$ that are orthogonal to $y$. Thus, if $\KK$
is the convex hull of finitely many points of $\mathbb Z ^d$, then
$$ \max_{x\in \KK} \la x, y \ra - \min_{x\in \KK} \la x, y \ra$$ 
is an integer namely, it is less by one than the number of lattice hyperplanes of $\mathbb Z ^d$ that intersect $\KK$ and are orthogonal to $y$. Now, we are ready to state the following conjecture of Hausel and the author (\cite{BeHa91}). 

\begin{con}\label{Bezdek-Hausel}
Let $\KK$ be a convex body in $\E^d$. Let $H_1, \dots, H_N$ be hyperplanes in $\E^d$ such that  
$$
   \KK \cap  \mathbb Z ^d \subset \bigcup _{i=1}^N H_i. 
$$
Then 
$$
   N\geq w(\KK, \mathbb Z ^d)-d.  
$$
 \end{con}

Properly translated copies of cross-polytopes, described in \cite{BeHa91}, show that if true, then the above inequality is best possible. 

The special case, when $N=0$, is of independent interest. Namely, it seems reasonable to conjecture (see also \cite{BaLiPaSz99}) that if $\KK$ is an integer point free convex body in $\E^d$, then $w(\KK, \mathbb Z ^d)\le d$. On the one hand, this has been proved by Banaszczyk \cite{Ba96} for ellipsoids. On the other hand, for general convex bodies containing no integer points, Banaszczyk, Litvak, Pajor and Szarek \cite{BaLiPaSz99} have proved the inequality $w(\KK, \mathbb Z ^d)\le C\cdot d^{\frac{3}{2}}$, where $C$ is an absolute positive constant. This improves an earlier result of Kannan and Lov\'asz \cite{KaLo88}.  

Although Conjecture~\ref{Bezdek-Hausel} is still open we have the following partial results published recently. Improving the estimates of \cite{BeHa91}, Talata \cite{Ta97} has succeeded to derive a proof of the following inequality. 

\begin{theorem}\label{Talata}
Let $\KK$ be a convex body in $\E^d$. Let $H_1, \dots, H_N$ be hyperplanes in $\E^d$ such that  
$$
   \KK \cap  \mathbb Z ^d \subset \bigcup _{i=1}^N H_i. 
$$
Then 
$$
   N\geq c\cdot\frac{w(\KK, \mathbb Z ^d)}{d}-d,  
$$
where $c$ is an absolute positive constant. 

\end{theorem}

In the paper \cite{BeLi09}, Litvak and the author have shown that the plank theorem of Ball \cite{B91} implies a slight improvement on the above inequality for centrally symmetric convex bodies whose lattice width is at most quadratic in dimension. (Actually, this approach is different from Talata's technique and can lead to a somewhat even stronger inequality in terms of the relevant basic measure of the given convex body. For more details on this we refer the interested reader to \cite{BeLi09}.)

\begin{theorem}\label{Bezdek-Litvak-2}
Let $\KK$ be a centrally symmetric convex body in $\E^d$. Let $H_1,$ $\dots, H_N$ be hyperplanes in $\E^d$ such that  
$$
   \KK \cap  \mathbb Z ^d \subset \bigcup _{i=1}^N H_i. 
$$
Then 
$$
   N\geq c\cdot \frac{w(\KK, \mathbb Z ^d)}{d\ln(d+1)},  
$$
where $c$ is an absolute positive constant. 
\end{theorem}
 
 Motivated by Conjecture~\ref{Bezdek-Hausel} and by a conjecture of Corzatt \cite{Co85} (according to which if in the plane the integer points of a convex domain can be covered by $N$ lines, then those integer points can also be covered by $N$ lines having at most four different slopes) Brass \cite{BrMoPa05} has raised the following related question.
 
\begin{prob}
For every positive integer $d$ find the smallest constant $c(d)$ such that if the integer points of a convex body in $\E^d$ can be covered by $N$ hyperplanes, then those integer points can also be covered by $c(d)\cdot N$ parallel hyperplanes.
\end{prob}

 Theorem~\ref{Talata} implies that $c(d)\le c\cdot d^2$ for convex bodies in general and for centrally symmetric convex bodies Theorem~\ref{Bezdek-Litvak-2} yields the somewhat better upper bound 
$c\cdot d\ln(d+1)$.

\section{Partial coverings by planks}

It seems that the following variant of Tarski's plank problem hasn't been yet considered: Let $\CC$ be a convex body of minimal width $w>0$ in $\E^d$. Moreover, let $w_1>0, w_2>0, \dots , w_n>0$ be given with $w_1+w_2+\dots+w_n<w$. Then find the arrangement of $n$ planks say, of $\PP_1, \PP_2,\dots , \PP_n$, of width $w_1, w_2, \dots , w_n$ in $\E^d$ such that their union covers the largest volume subset of $\CC$ that is for which  ${\rm vol}_d((\PP_1\cup\PP_2\cup\dots\cup\PP_n)\cap\CC )$ is as large as possible. As the following special case is the most striking form of the above problem, we are putting it forward as the main question of this section.

\begin{prob}
Let $\BB^d$ denote the unit ball centered at the origin ${\bf o}$ in $\E^d$. Moreover, let $w_1, w_2, \dots , w_n$ be positive real numbers satisfying the inequality $w_1+w_2+\dots+w_n<2$. Then prove or disprove that the union of the planks $\PP_1, \PP_2,\dots , \PP_n$ of width $w_1, w_2, \dots , w_n$ in $\E^d$ covers the largest volume subset of $\BB^d$ if and only if $\PP_1\cup\PP_2\cup\dots\cup\PP_n $ is a plank of width $w_1+w_2+\dots+w_n$ with ${\bf o}$ as a center of symmetry.
\end{prob}

Clearly, there is an affirmative answer to the above problem for $n=1$. Also, we have the following positive results. For the sake of completeness we include their short proofs.

\begin{theorem}\label{elso}
If $\PP_1$ and $\PP_2$ are planks in $\E^d, d\ge 2$ of width $w_1$ and $w_2$ having $0<w_1+w_2<2$, then $\PP_1\cup\PP_2$ covers the largest volume subset of $\BB^d$ if and only if $\PP_1\cup\PP_2$ is a plank of width $w_1+w_2$ possessing ${\bf o}$ as a center of symmetry.
\end{theorem}

\proof The following is an outline of a quite elementary proof. First, let us consider the case when $\PP_1$ and $\PP_2$ are planks in $\E^2$ of width $w_1$ and $w_2$ having $0<w_1+w_2<2$. We say, that $(\PP_1\cup\PP_2)\cap\BB^2$ is a crossing subset of $\BB^2$, if $\BB^2\setminus(\PP_1\cup\PP_2)$ consists of $4$ connected components. Now, it is not hard to see that among the crossing subsets (resp., non-crossing subsets) the only critical configuration with respect to maximizing the area is the one with $\PP_1$ and $\PP_2$ being perpendicular to each other and having ${\bf o}$ as a center of symmetry (resp., the one with $\PP_1\cup\PP_2$ being a plank of width $w_1+w_2$ and having ${\bf o}$ as a center of symmetry). Second, it is easy to check that between the two critical configurations the non-crossing one possesses a larger area, finishing the proof Theorem~\ref{elso} for $d=2$. Finally, if $\PP_1$ and $\PP_2$ are planks in $\E^d, d\ge 3$ of width $w_1$ and $w_2$ having $0<w_1+w_2<2$, then an application of the $2$-dimensional case of Theorem~\ref{elso}, just proved, to the $2$-dimensional flats of $\E^d$ that are parallel to the normal vectors of $\PP_1$ and $\PP_2$ followed by an integration of the areas of the corresponding sets sitting on the $2$-flats in question, yield the desired claim.  \kkk

\begin{theorem}\label{masodik}
Let $w_1, w_2, \dots , w_n$ be positive real numbers satisfying the inequality $w_1+w_2+\dots+w_n<2$. Then the union of the planks $\PP_1, \PP_2,\dots , \PP_n$ of width $w_1, w_2, \dots , w_n$ in $\E^3$ covers the largest volume subset of $\BB^3$ if and only if $\PP_1\cup\PP_2\cup\dots\cup\PP_n $ is a plank of width $w_1+w_2+\dots+w_n$ with ${\bf o}$ as a center of symmetry.
\end{theorem}

\proof Let $\PP_1, \PP_2,\dots , \PP_n$ be an arbitrary family of planks of width $w_1,$ $w_2, \dots , w_n$ in $\E^3$ and let $\PP$ be a plank of width $w_1+w_2+\dots+w_n$ with ${\bf o}$ as a center of symmetry. Moreover, let $S(x)$ denote the sphere of radius $x$ centered at ${\bf o}$. Now, recall the well-known fact that if $\PP(y)$ is a plank of width $y$ whose each boundary plane intersects $S(x)$, then ${\rm sarea} (S(x)\cap \PP(y) )=2\pi xy$, where ${\rm sarea}(\ .\ )$ refers to the surface area measure on $S(x)$. This implies in a straightforward way that
$${\rm sarea}[(\PP_1\cup\PP_2\cup\dots\cup\PP_n)\cap S(x) ]\le {\rm sarea}(\PP\cap S(x)),$$
and so,
$${\rm vol}_3((\PP_1\cup\PP_2\cup\dots\cup\PP_n)\cap\BB^3 )=\int_0^1 {\rm sarea}[(\PP_1\cup\PP_2\cup\dots\cup\PP_n)\cap S(x) ]\ dx \le $$
$$ \int_0^1 {\rm sarea}(\PP\cap S(x))\ dx ={\rm vol}_3(\PP\cap\BB^3 ),$$ 
finishing the proof of the "if" part of Theorem~\ref{masodik}. Actually, a closer look of the above argument gives a proof of the "only if" part as well. \kkk

As an immediate corollary we get the following statement.

\begin{cor}\label{harmadik}
If $\PP_1, \PP_2$ and $\PP_3$ are planks in $\E^d, d\ge 3$ of width $w_1, w_2$ and $w_3$ satisfying $0<w_1+w_2+w_3<2$, then $\PP_1\cup\PP_2\cup\PP_3$ covers the largest volume subset of $\BB^d$ if and only if $\PP_1\cup\PP_2\cup\PP_3$ is a plank of width $w_1+w_2+w_3$ having ${\bf o}$ as a center of symmetry.
\end{cor}

\proof Indeed an application of Theorem~\ref{masodik} to the $3$-dimensional flats of $\E^d$ that are parallel to the normal vectors of $\PP_1, \PP_2$ and $\PP_3$ followed by an integration of the volumes of the corresponding sets lying in the $3$-flats in question, yield the desired claim.  \kkk

In general, we have the following estimate that one can derive from Bang's paper \cite{Ba51} as follows. In order to state it properly we introduce two definitions.

\begin{definition}
{\rm Let $\CC$ be a convex body in $\E^d$ and let $m$ be a positive integer. Then let ${\cal T}_{\CC, d}^{m}$ denote the family of all sets in $\E^d$ that can be obtained as the intersection of at most $m$ translates of $\CC$ in $\E^d$.}
\end{definition}

\begin{definition}
{\rm Let $\CC$ be a convex body of minimal width $w>0$ in $\E^d$ and let $0<x\le w$ be given. Then for any non-negative integer $n$ let
$$v_d(\CC, x, n ):=\min \{ {\rm vol}_d(\QQ)\ |\ \QQ\in {\cal T}_{\CC, d}^{2^n}\ {\rm and}\ w(\QQ)\ge x \ \}.$$}
\end{definition}

Now, we are ready to state the theorem which although was not published by Bang in \cite{Ba51}, it follows from his proof of Tarski's plank conjecture.

\begin{theorem}\label{Bang} 
Let $\CC$ be a convex body of minimal width $w>0$ in $\E^d$. Moreover, let $\PP_1, \PP_2,\dots , \PP_n$ be planks of width $w_1, w_2, \dots , w_n$ in $\E^d$ with $w_0=w_1+w_2+\dots+w_n<w$. Then
$${\rm vol}_d(\CC\setminus(\PP_1\cup\PP_2\cup\dots\cup\PP_n))\ge v_d(\CC, w-w_0, n),$$
that is
$${\rm vol}_d((\PP_1\cup\PP_2\cup\dots\cup\PP_n)\cap\CC )\le {\rm vol}_d(\CC )-v_d(\CC, w-w_0, n).$$
\end{theorem}

Clearly, the first inequality above implies (via an indirect argument) that if the planks $\PP_1, \PP_2,\dots , \PP_n$ of width $w_1, w_2, \dots , w_n$ cover the convex body $\CC$ in $\E^d$, then $w_1+w_2+\dots+w_n\ge w$. Also, as an additional observation we mention the following statement, that on the one hand, can be derived from Theorem~\ref{Bang} in a straightforward way, on the other hand, represents the only case when the estimate in Theorem~\ref{Bang} is sharp.

\begin{cor}
Let $\TT$ be an arbitrary triangle of minimal width (i.e. of minimal height) $w>0$ in $\E^2$. Moreover, let $w_1, w_2, \dots , w_n$ be positive real numbers satisfying the inequality $w_1+w_2+\dots+w_n<w$. Then the union of the planks $\PP_1, \PP_2,\dots , \PP_n$ of width $w_1, w_2, \dots , w_n$ in $\E^2$ covers the largest area subset of $\TT$ if $\PP_1\cup\PP_2\cup\dots\cup\PP_n $ is a plank of width $w_1+w_2+\dots+w_n$ sitting on the side of $\TT$ with height $w$.
\end{cor}

\medskip

\section{Linking the Kakeya-P\'al and the Blasch\-ke-Lebesgue problems to the partial covering problem}

Recall that the {\it Kakeya-P\'al problem} is about minimizing the volume of convex bodies of given minimal width $w>0$ in $\E^d$. For short reference let $\KK_{{KP}}^{w, d}$ denote any of the minimal volume convex bodies in the Kakeya-P\'al problem. (Actually, Kakeya phrased his question in 1917 as follows: what is the smallest area of a convex set within which one can rotate a needle by $180^{{\rm o}}$.) P\'al \cite{P21} has solved this problem for $d=2$ by showing that the smallest area convex domain of minimal width $w>0$ is a regular triangle of height $w$. As it is well-known, the Kakeya-P\'al problem is unsolved in higher dimensions (for more details on this see for example \cite{CCG96}). Thus, the following is an immediate corollary of Theorem~\ref{Bang}.

\begin{cor}\label{upper-bound} 
Let $\CC$ be a convex body of minimal width $w>0$ in $\E^d$. Moreover, let $\PP_1, \PP_2,\dots , \PP_n$ be planks of width $w_1, w_2, \dots , w_n$ in $\E^d$ with $w_0=w_1+w_2+\dots+w_n<w$. Then
$${\rm vol}_d((\PP_1\cup\PP_2\cup\dots\cup\PP_n)\cap\CC )\le {\rm vol}_d(\CC )-{\rm vol}_d(\KK_{{KP}}^{w-w_0, d}).$$
\end{cor}

It seems that the best lower bound for the Kakeya-P\'al problem is due to Firey \cite{Fi65} claiming that ${\rm vol}_d(\KK_{{KP}}^{w, d})\ge f(d)w^d$ with $f(d)=\frac{2}{\sqrt{3}\cdot d!}$. Corollary~\ref{upper-bound} suggests to further investigate and improve Firey's inequality  for $d\ge 3$. (For $d=2$ the inequality in question is identical to P\'al's result \cite{P21} and so, it is optimal.) Here, we claim the following improvement.

\begin{theorem}\label{first}
Let $\CC$ be a convex body of minimal width $w>0$ in $\E^d$. Moreover, for each odd integer $d\ge 3$ let $g(d)=\sqrt{\frac{3\cdot\pi^{d-3}\cdot(d+1)!!}{2^{d-2}\cdot (d!!)^5} }$ and for each even integer $d\ge 4$ let
$g(d)=\sqrt{\frac{3\cdot \pi^{d-3}\cdot (d+2)!! }{(d+1)^2\cdot(d!!)^2\cdot ((d-1)!!)^3 } }$.

Then
$$ {\rm vol}_d(\CC )\ge g(d)w^d>f(d)w^d$$
for all $d\ge 3$.
\end{theorem}

\proof We outline the proof by describing its main idea and by leaving out the more or less straightforward but, somewhat lengthy computations. First, we need the following result of Steinhagen \cite{St22}. Let $\CC$ be a convex body of minimal width $w>0$ in $\E^d$. Moreover, for each odd integer $d\ge 3$ let $r(d)=\frac{1}{2\sqrt{d}}$
and for each even integer $d\ge 2$ let $r(d)=\frac{\sqrt{d+2}}{2(d+1)}$. Then the inradius $r$ of $\CC$ (which is the radius of the largest ball lying in $\CC$) is always at least as large as $r(d)w$. Second, recall Kubota's formula \cite{Bo87} according to which  
$${\rm svol}_{d-1}( {\rm bd}(\CC))= \frac{1}{{\rm vol}_{d-1}(\BB^{d-1})}\int_{\S^{d-1}}{\rm vol}_{d-1}(\CC \big| {\bf x})\ d{\bf x},$$  
where ${\rm bd}(\ .\ )$ (resp., ${\rm svol}_{d-1}(\ .\ )$) stands for the boundary (resp., $(d-1)$-dimensional surface volume) of the corresponding set and $\S^{d-1}= {\rm bd}(\BB^{d})$ moreover, $\CC \big| {\bf x}$ denotes the orthogonal projection of $\CC$ onto the hyperplane passing through ${\bf o}$ with normal vector ${\bf x}$. Thus, Steinhagen's theorem and Kubota's formula imply in a straighforward way that
$$ {\rm vol}_d(\CC )\ge  \frac{r(d)w}{d}  {\rm svol}_{d-1}( {\rm bd}(\CC))\ge 
\frac{r(d)w {\rm vol}_{d}(\BB^{d})}{{\rm vol}_{d-1}(\BB^{d-1})}\min_{{\bf x}\in \S^{d-1} }\{{\rm vol}_{d-1}(\CC \big| {\bf x})\}  .$$ 
Finally, as $\CC \big| {\bf x}$ is a $(d-1)$-dimensional convex body of minimal width at least $w$ for all $ {\bf x}\in \S^{d-1} $, therefore the above inequality, repeated in a recursive way for lower dimensions, leads to the desired inequality claimed in Theorem~\ref{first}. \kkk

\begin{remark}
{\rm For comparison we mention that $g(3)=\frac{2}{9}=0.2222 \dots >f(3)=\frac{1}{3\sqrt{3}}=0.1924\dots $ (resp., $g(4)=\sqrt{\frac{2\pi}{75}}=0.2894\dots >f(4)=\frac{1}{12\sqrt{3}}=0.0481\dots $). Also, recall that Heil \cite{He78} has constructed a convex body in $\E^3$ of minimal width $1$ and of volume $0.298...$.}
\end{remark}

 Corollary~\ref{upper-bound} can be further improved when $\CC$ is a Euclidean ball. The details are as follows.

First, recall that the {\it Blaschke-Lebesgue problem} is about finding the minimum volume convex body of constant width $w>0$ in $\E^d$. In particular, the Blaschke-Lebesgue theorem states that among all convex domains of constant width $w$, the Reuleaux triangle of width $w$ has the smallest area, namely $\frac{1}{2}(\pi-\sqrt{3})w^2$. W. Blaschke \cite{Bla15} and H. Lebesgue \cite{Leb14} were the first to show this and the succeding decades have seen other works published on different proofs of that theorem. For a most recent new proof, and for a survey on the state of the art of different proofs of the Blaschke-Lebesgue theorem, see the elegant paper of E. M. Harrell \cite{Ha02}. Here we note that the Blaschke-Lebesgue problem is unsolved in three and more dimensions. Even finding the $3$-dimensional set of least volume presents
formidable difficulties. On the one hand, Chakerian \cite{Ch66} proved that any convex body of constant width $1$ in $\E^3$ has volume at least
$\frac{\pi (3\sqrt{6}-7) }{3}=0.365...$. On the other hand, it has been conjectured by Bonnesen and Fenchel \cite{Bo87} that Meissner's $3$-dimensional generalizations of the Reuleaux triangle of volume $\pi (\frac{2}{3}-\frac{1}{4}\sqrt{3}\arccos (\frac{1}{3}))=0.420...$ are the only extramal sets in $\E^3$.

For our purposes it will be useful to introduce the notation ${\KK}^{w, d}_{{BL}}$ (resp., ${\overline{\KK}}^{w, d}_{{BL}}$)
for a convex body of constant width $w$ in $\E^d$ having minimum volume (resp., surface volume). One may call ${\KK}^{w, d}_{{BL}}$ (resp., ${\overline{\KK}}^{w, d}_{{BL}}$) a Blaschke-Lebesgue-type convex body with respect to volume (resp., surface volume). Note that for $d=2, 3$ one may choose ${\KK}^{w, d}_{{BL}}={\overline{\KK}}^{w, d}_{{BL}}$ however, this is likely not to happen for $d\ge 4$. (For more details on this see \cite{Ch66}.) As an important note we mention that Schramm \cite{Sc88} has proved the inequality
$${\rm vol}_d( {\KK}^{w, d}_{{BL}})\ge \bigg(\sqrt{3+\frac{2}{d+1}}-1 \bigg)^d(\frac{w}{2})^d {\rm vol}_d(\BB^d ),$$ which gives the best lower bound for all $d>4$. By observing that the orthogonal projection of a convex body of constant width $w$ in $\E^d$ onto any hyperplane of $\E^d$ is a $(d-1)$-dimensional convex body of constant width $w$ one obtains from the previous inequality of Schramm the following one: 
$${\rm svol}_{d-1}( {\rm bd}({\overline{\KK}}^{w, d}_{{BL}}))\ge d\bigg(\sqrt{3+\frac{2}{d}}-1 \bigg)^{d-1}(\frac{w}{2})^{d-1} {\rm vol}_{d-1}(\BB^{d-1} ){\rm vol}_d(\BB^d ).$$ 

Second, let us recall that if $X$ is finite (point) set lying in the interior of a unit ball in $\E^d$, then the intersection of the (closed) unit balls of $\E^d$ centered at the points of $X$ is called a ball-polyhedron and it is denoted by $\BB[X]$. (For an extensive list of properties of ball-polyhedra see the recent paper \cite{Be07}.) Of course, $\BB[X]$ makes sense to introduce also for sets $X$ that are not finite but, in those cases we get sets that are typically not ball-polyhedra. 

Now, we are ready to state our theorem.

\begin{theorem}\label{second}
Let $\BB[X]\subset \E^d$ be a ball-polyhedron of minimal width $x$ with $1\le x<2$. Then
$${\rm vol}_d(\BB[X])\ge {\rm vol}_d({\KK}^{2-x, d}_{{BL}})+{\rm svol}_{d-1}({\rm bd}({\overline{\KK}}^{2-x, d}_{{BL}}))\cdot (x-1)+{\rm vol}_d(\BB^d )\cdot (x-1)^d.$$
\end{theorem}

\proof Recall that if if $X$ is finite set lying in the interior of a unit ball in $\E^d$, then we can talk about its spindle convex hull ${\rm conv}_{\rm s}(X)$, which is simply the intersection of all (closed) unit balls of $\E^d$ that contain $X$ (for more details see \cite{Be07}). The following statement can be obtained by combining Corollary 3.4 of \cite{Be07} and Proposition 1 of \cite{Be06}.

\begin{lemma}\label{basic}
Let $X$ be a finite set lying in the interior of a unit ball in $\E^d$. Then

\item (i) ${\rm conv}_{\rm s}(X)=\BB \big[\BB[X]\big]$ and therefore $\BB[X]=\BB\big[{\rm conv}_{\rm s}(X)\big]$;

\item (ii) the Minkowski sum $\BB[X]+ {\rm conv}_{\rm s}(X)$ is a convex body of constant width $2$ in $\E^d$ and so, $w(\BB[X])+{\rm diam}\big({\rm conv}_{\rm s}(X) \big)=2$, where ${\rm diam}(\ .\ )$ stands for the diameter of the corresponding set in $\E^d$.

\end{lemma}

By part $(ii)$ of Lemma~\ref{basic} ${\rm diam}\big({\rm conv}_{\rm s}(X) \big)\le 2-x$. This implies, via a classical theorem of convexity (see for example \cite{Bo87}), the existence of a convex body $\LL$ of constant width $(2-x)$ in $\E^d$ with ${\rm conv}_{\rm s}(X)\subset \LL$. Hence, using part $(i)$ of Lemma~\ref{basic}, we get that $\BB[{\LL} ]\subset{\BB}[X]={\BB}\big[{\rm conv}_{\rm s}(X)\big]$. Finally, notice that as $\LL$ is a convex body of constant width $(2-x)$ therefore $\BB[{\LL}]$
is in fact, the outer-parallel domain of $\LL $ having radius $(x-1)$ (that is $\BB[{\LL}]$ is the union of all $d$-dimensional (closed) balls of radii $(x-1)$ in $\E^d$ that are centered at the points of $\LL $). Thus, 
$${\rm vol}_d(\BB[X])\ge {\rm vol}_d\big(\BB[{\LL}]\big)=
{\rm vol}_d(\LL )+{\rm svol}_{d-1}({\rm bd}(\LL ))\cdot (x-1)+{\rm vol}_d(\BB^d )\cdot (x-1)^d.$$

The inequality above together with the following obvious ones  
$${\rm vol}_d(\LL )\ge {\rm vol}_d({\KK}^{2-x}_{{BL}_d}) {\rm \ and \ } {\rm svol}_{d-1}({\rm bd}(\LL ))\ge {\rm svol}_{d-1}({\rm bd}({\overline{\KK}}^{2-x}_{{BL}_d}))$$ 
imply Theorem~\ref{second} in a straightforward way. \kkk

Thus, Theorem~\ref{Bang} and Theorem~\ref{second} imply the following immediate estimate.

\begin{cor}
Let $\BB^d$ denote the unit ball centered at the origin ${\bf o}$ in $\E^d, d\ge 2$. Moreover, let $\PP_1, \PP_2,\dots , \PP_n$ be planks of width $w_1, w_2, \dots , w_n$ in $\E^d$ with $w_0=w_1+w_2+\dots+w_n\le 1$. Then
$${\rm vol}_d((\PP_1\cup\PP_2\cup\dots\cup\PP_n)\cap\BB^d )\le {\rm vol}_d(\BB^d )-v_d(\BB^d, 2-w_0, n)\le $$ 
$$\big(1-(1-w_0)^d\big){\rm vol}_d(\BB^d )- {\rm vol}_d({\KK}^{w_0, d}_{{BL}})-{\rm svol}_{d-1}({\rm bd}({\overline{\KK}}^{w_0, d}_{{BL}}))\cdot \big(1-w_0)$$
\end{cor}

\medskip

\section{Strengthening the plank theorems of Ball and Bang}

Recall that Ball (\cite{B91}) generalized the plank theorem of Bang (\cite{Ba50}, \cite{Ba51}) for coverings of balls by planks in Banach spaces (where planks are defined with the help of linear functionals instead of inner product). This theorem was further strengthened by Kadets \cite{Ka05} for Hilbert spaces as follows. Let $\CC$ be a closed convex subset with non-empty interior in the real Hilbert space $\H$ (finite or infinite dimensional). We call $\CC$ a {\it convex body} of $\H$. Then let ${\rm r}(\CC)$ denote the supremum of the radii of the balls contained in $\CC$. (One may call ${\rm r}(\CC)$ the {\it inradius} of $\CC$.) Planks and their widths in $\H$ are defined with the help of the inner product of $\H$ in the usual way. Thus, if $\CC$ is a convex body in $\H$ and $\PP$ is a plank of $\H$, then the width ${\rm w}(\PP)$ of $\PP$ is always at least as large as $2\cdot {\rm r}(\CC\cap\PP)$. Now, the main result of \cite{Ka05} is the following. 

\begin{theorem}\label{Kadets-theorem}
Let the ball $\BB$ of the real Hilbert space $\H$ be covered by the convex bodies $\CC_1, \CC_2, \dots , \CC_n$ in $\H$. Then
$$\sum_{i=1}^n {\rm r}(\CC_i)\ge {\rm r}(\BB).$$
\end{theorem}

We note that an independent proof of the $2$-dimensional Euclidean case of Theorem~\ref{Kadets-theorem} can be found in \cite{Be-A07}. Kadets (\cite{Ka05}) proposes to investigate the analogue of Theorem~\ref{Kadets-theorem} in Banach spaces. Thus, an affirmative answer to the following problem would improve the plank theorem of Ball.

\begin{prob}
Let the ball $\BB$ be covered by the convex bodies $\CC_1, \CC_2, \dots , \CC_n$ in an arbitrary Banach space. Prove or disprove that
$$\sum_{i=1}^n {\rm r}(\CC_i)\ge {\rm r}(\BB).$$
\end{prob}

In \cite{Be-A03}, A. Bezdek proposes another way for improving the plank theorem of Bang.

\begin{con}\label{A.Bezdek-theorem}
For each convex region $\CC$ in $\E^2$ there exists an $\epsilon>0$ such that if $\epsilon\CC$ lies in the interior of $\CC$ and the annulus $\CC\setminus\epsilon\CC$ is covered by finitely many planks, then the sum of the widths of the planks is at least the minimal width of $\CC$.
\end{con}

The following theorem, proved in \cite{Be-A03}, supports Conjecture~\ref{A.Bezdek-theorem}.

\begin{theorem}
Let $\CC\subset\E^2$ be a unit square and let $\epsilon =1-\frac{1}{\sqrt{2}}=0.292...$. If $\epsilon\CC$ lies in the interior of $\CC$ and the annulus $\CC\setminus\epsilon\CC$ is covered by finitely many planks, then the sum of the widths of the planks is at least $1$. 
\end{theorem}

\vspace{1cm}

\medskip

\noindent
K\'aroly Bezdek,
Department of Mathematics and Statistics,
2500 University drive N.W.,
University of Calgary, Calgary, Alberta, Canada, T2N 1N4.
\newline
{\sf e-mail: bezdek@math.ucalgary.ca}

\end{document}